\documentclass[final]{siamltex}
\usepackage{hyperref}
\hypersetup{%
    pdftitle={A multigrid method for CutFEM and its implementation on GPU},
    pdfauthor={Cu Cui, Guido Kanschat},
    pdfsubject={65N55; 65Y05; 65Y10},
    pdfkeywords={unfitted finite elements; multigrid; vertex-patch smoother; GPU implementation}
    }
\usepackage{times}
\usepackage{graphicx} 

\usepackage{amsmath}
\usepackage{amsfonts}

\usepackage{tikz,tikzscale}
\usepackage{pgfplotstable}
\pgfplotsset{compat=1.18}
\usetikzlibrary{calc}
\usetikzlibrary{shadings}
\usetikzlibrary{shapes.geometric}
\usetikzlibrary{shadows.blur}
\usetikzlibrary{positioning}
\usetikzlibrary{backgrounds}
\usetikzlibrary{matrix}
\usetikzlibrary{intersections}

\usepackage[notref,notcite]{showkeys}
\usepackage[mode=multiuser]{fixme}
\fxusetheme{color}
\FXRegisterAuthor{cc}{envfb}{CC}
\FXRegisterAuthor{gk}{envgk}{GK}

\usepackage{multirow}
\usepackage{booktabs}
\usepackage{stmaryrd}

\usepgfplotslibrary{fillbetween}
\usepgfplotslibrary{colormaps}
\makeatletter
\definecolor{gnuplot@orange}{RGB}{229,158,0}
\definecolor{gnuplot@purple}{RGB}{148,0,212}
\definecolor{gnuplot@lightblue}{RGB}{87,181,232}
\definecolor{gnuplot@green}{RGB}{0,158,75}
\definecolor{gnuplot@darkblue}{RGB}{0,115,179}
\definecolor{gnuplot@yellow}{RGB}{240,227,66}
\pgfplotscreateplotcyclelist{colorGPL}{%
gnuplot@darkblue,every mark/.append style={fill=gnuplot@darkblue!80!black},mark=o\\%
red!80!white,every mark/.append style={fill=red!80!black},mark=square\\%
gnuplot@green,every mark/.append style={fill=gnuplot@green!50!black},mark=otimes*\\%
gnuplot@orange,every mark/.append style={fill=gnuplot@orange!80!black,mark size=2.5pt},mark=x\\%
gnuplot@purple,mark=diamond\\%
black,densely dashed,every mark/.append style={solid,fill=gnuplot@darkblue!80!black},mark=square*\\%
gnuplot@lightblue,densely dashed,every mark/.append style={solid,fill=gnuplot@lightblue!30!black},mark=triangle*\\%
red!80!white,densely dashed,every mark/.append style={solid,fill=red!40!white},mark=oplus*\\%
}
\makeatother

\def\mesh{\mathbb M}
\def\cell{T}
\def\vertices{\mathbb V}
\def\vertex{X}
\def\Q{\mathbb Q}
\def\prol#1{I^\uparrow_{#1}}
\def\rest#1{I^\downarrow_{#1}}
\def\vA{\mathbf{A}}
\def\vQ{\mathbf{Q}}
\def\vb{\mathbf{b}}
\def\vx{\mathbf{x}}
\def\vy{\mathbf{y}}

\title{A multigrid method for CutFEM and its implementation on GPU}
\author{Cu Cui\footnotemark[2] \and Guido Kanschat\footnotemark[2]}

\begin{document}

\maketitle

\renewcommand{\thefootnote}{\fnsymbol{footnote}}

\footnotetext[2]{Interdisciplinary Center for Scientific Computing (IWR), Faculty for Engineering Sciences, Heidelberg University, Im Neuenheimer Feld 205, 69120 Heidelberg, Germany, \texttt{cu.cui@alumni.uni-heidelberg.de}, \texttt{kanschat@uni-heidelberg.de}.\\
Cu Cui was supported by the China Scholarship Council (CSC) under grant NO. 202106380059.
}

\begin{abstract}
    We present a multigrid method for an unfitted finite element discretization of the Dirichlet boundary value problem. The discretization employs Nitsche's method to implement the boundary condition and additional face based ghost penalties for stabilization. We apply standard intergrid operators, relying on the fact that the relevant domain of computation does not grow under mesh refinement. The smoother is a parallel implementation of the multiplicative vertex-patch smoother with inconsistent treatment of ghost penalties. Our computational results show that we obtain a fast converging method. Furthermore, runtime comparison to fitted methods show that the losses are moderate although many optimizations for Cartesian vertex patches cannot be applied on cut patches.
\end{abstract}



\section{Introduction}

Our goal in this article is combining the advantages of unfitted finite element discretizations with state-of-the-art GPU implementations of finite element solvers.
Unfitted finite element methods provide a powerful framework for handling complex geometries as well as moving boundaries and interfaces while enabling the use of computationally efficient tensor techniques in the bulk of the domain. In this work, we establish a geometric multigrid technique based on vertex-patch smoothers to obtain fast solvers and efficient preconditioners for such discretizations.

We focus on the realization of unfitted finite elements by the CutFEM method, going back to the article by Hansbo and Hansbo~\cite{HansboHansbo02}. For the handling of instabilities due to small cuts, we rely on ``ghost penalties'' introduced and analyzed in~\cite{burman2010ghost,BurmanHansbo12}.

A multigrid method for the unfitted interface problem was presented in~\cite{LudescherGrossReusken20}. There, a Gauss-Seidel smoother is used inside the two subdomains and combined with a global solver for the degrees of freedom coupling over the interface.
They obtain good iteration counts independent of contrast, such that their method can be expected to work for the unfitted boundary value problem as well.
Nevertheless, such a global solver for the interface degrees of freedom necessitates data structures which impede the flow of information on the GPU, and hence we prefer to avoid it.
In the thesis~\cite{Ludescher20}, a basis change at the interface is presented which renders good convergence rates without a global solution involving the interface matrix.
Additionally, results for a $p$-multigrid method for higher order discretizations are presented in the thesis, employing face-patch solvers.
Furthermore, in~\cite{GrossReusken23}, a domain decomposition preconditioner for the CutFEM interface problem is presented and analyzed. Albeit their estimates are not robust, they also present computational results for the Dirichlet boundary value problem. Another method related to our work was developed in~\cite{PrenterVerhooselBrummelen19} for unfitted isogeometric methods. There, a cell-based additive Schwarz method for all cut cells is combined with a Jacobi method on the uncut cells. Since they do not combine this preconditioner with a coarse grid, they exhibit the usual dependence of the condition number for second order problems, namely $\kappa \sim h^{-2}$, but independent of the size and location of the cuts.

The vertex-patch smoother as a particular implementation of a multilevel domain decomposition method was introduced by Pavarino in~\cite{Pavarino94additive}.
In its additive and multiplicative form, it has proven to be a robust smoother for many applications~\cite{ArnoldFalkWinther00,KanschatLazarovMao17,MeggendorferKanschatKraus24}, can be implemented in an efficient way~\cite{BrubeckFarrell22,WitteArndtKanschat21}, and provides the computational intensity needed to exploit modern computers~\cite{WichrowskiMunchKronbichlerKanschat25}. In a series of publications, we have established that it enables implementations on GPU with high utilization~\cite{Cui24,CuiGrosseBleyKanschatStrzodka25,CuiKanschat24IP,CuiKanschat25Stokes}, exemplified with the CUDA and Tensor cores of the Nvidia A100, respectively.

Our numerical experiments demonstrate that the proposed approach not only achieves robustness in GPU implementations for unfitted methods but also delivers performance comparable to that observed on standard Cartesian geometries. This advancement represents a significant step toward extending GPU-accelerated solvers to complex geometries and unfitted methods, underscoring the potential for combining multigrid techniques with Schwarz smoothers to efficiently manage the computational demands of CutFEM on modern hardware.

\section{Model Problem and Multigrid Solver}

\begin{figure}[tp]
\centering
\begin{tikzpicture}[scale=0.6]
    \fill [fill=blue!25!white] (2,2) rectangle (6,6);
    \fill [fill=blue!25!white] (2,1) rectangle (6,2);
    \fill [fill=blue!25!white] (2,6) rectangle (6,7);
    \fill [fill=blue!25!white] (1,2) rectangle (2,6);
    \fill [fill=blue!25!white] (6,2) rectangle (7,6);
    \draw[color=blue,very thick] (4,4) circle (2.4);
    \draw[step=1cm, black, thin] (0,0) grid (8,8);

\end{tikzpicture}
\caption{Circular domain $\Omega$ on top of a uniform, unfitted background mesh $\mesh_\ell$, and shaded the domain $\Omega_\ell$ consisting of all cells having a nonempty intersection with $\Omega$.}\label{fig:cutfem}
\end{figure}
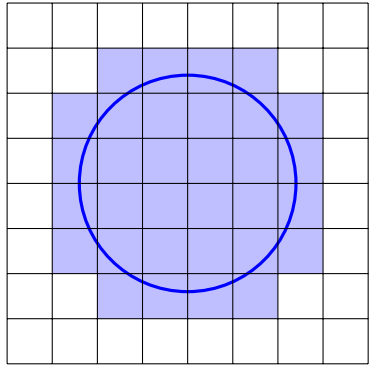

As a model problem, we consider the Poisson equation on the domain $\Omega$ with a uniform, unfitted background mesh, as illustrated in Figure~\ref{fig:cutfem}.

\subsection{Discretization}
The problem is discretized using the continuous Galerkin (CG) method, and the Dirichlet boundary condition is enforced weakly via the Nitsche method~\cite{HansboHansbo02,nitsche1971variationsprinzip}.
To this end, we first introduce a hierarchy of Cartesian meshes
\begin{gather*}
    \mesh_0 \sqsubset \mesh_1 \sqsubset\dots\sqsubset \mesh_L,
\end{gather*}
where the symbol ``$\sqsubset$'' indicates nestedness, that is, the cells of the mesh on the left are unions of cells of the mesh on the right. For each level $\ell$, we introduce the subsets $\mesh_{\ell,\Omega}$ as the set of all cells of $\mesh_\ell$ with nonempty intersection with the domain $\Omega$ and $\mesh_{\ell,\Gamma}$ as those intersected by the boundary $\Gamma = \partial\Omega$.

For each level $0\le\ell\le L$, we introduce the finite element space
\begin{gather}
    \label{eq:spaces}
    V_\ell = V_{\ell,p}
    = \left\{ v \in H^1(\Omega) \middle|
    v_{|\cell} \in \Q_p(\cell) \; \forall \cell\in\mesh_{\ell,\Omega}
    \right\}.
\end{gather}
Here, $\Q_p(\cell)$ is the space of tensor product polynomials of degree up to $p$ with a nodal basis defined by interpolation in Gauss-Lobatto points. By concatenation, we obtain a nodal basis $\{\phi_{\ell,i}\}$ for the space $V_\ell$. Note that the support of the space $V_\ell$ is larger than $\Omega$. As the polynomials are defined on whole cells, it is the union of all cells in $\mesh_{\ell,\Omega}$. We call this set $\Omega_\ell$.
Following~\cite{HansboHansbo02}, we introduce on $V_\ell$ the bilinear form
\begin{equation}\label{eq:cutfem_nitsche}
  a_\ell(u_\ell,v_\ell)  =  (\nabla u_\ell, \nabla v_\ell)_\Omega
                  - (\partial_n u_\ell, v_\ell)_\Gamma
                  - (u_\ell, \partial_n v_\ell)_\Gamma
                  + \left(\frac{\gamma_D}{h_\ell} u_\ell, v_\ell \right)_\Gamma,
\end{equation}
where \(\partial_n = n \cdot \nabla\) denotes the normal derivative in the direction of the unit outward normal \(n\). Furthermore, \(h_\ell\) denotes the characteristic element size and \(\gamma_D\) is the Nitsche penalty parameter, which here may depend on the particular cut cell. Note that the integration over \(\Omega\) involves only those cells which are completely inside the domain and parts of the cells which overlap the domain boundary.

It has been established, see~\cite{burman2010ghost}, that the weak formulation~\eqref{eq:cutfem_nitsche} exhibits numerical instability due to the small cut cell problem. Two primary challenges arise: (i) the Nitsche penalty parameter \(\gamma_D\) used to enforce the boundary condition would have to be arbitrarily large, see also~\cite{de2018note}, adversely affecting stability, and (ii) the stiffness matrices can become arbitrarily ill-conditioned when the shape functions have only a small support on cut cells.

To address these issues, stabilization techniques have been proposed in the literature. In this work, we adopt a face-based ghost penalty method~\cite{BurmanHansbo12} to enhance stability and improve conditioning. The ghost penalty stabilization augments the original bilinear form with an additional, consistent term:
\begin{gather}
  \label{eq:cutfem-ghost}
  A_\ell(u_\ell,v_\ell)  = a_\ell(u_\ell,v_\ell) + g_\ell(u_\ell, v_\ell),
\end{gather}
The stabilization term \(g_\ell(u_\ell, v_\ell)\) is defined on all faces between cut cells and all other cells in $\mesh_{\ell,\Omega}$. More precisely, let for two neighboring cells $F(\cell_1,\cell_2)$ be their joint face of codimension one. Then, define the set of ghost faces
\begin{gather*}
    \mathbb F_G = \bigl\{ F(\cell_1,\cell_2)\bigm|
    \cell_1,\cell_2\in\mesh_{\ell,\Omega},\;
    \cell_1\in\mesh_{\ell,\Gamma}\lor\cell_2\in\mesh_{\ell,\Gamma}
    \bigr\},
\end{gather*}
and the ghost penalty
\begin{equation*}
  g_\ell(u_\ell, v_\ell) = 
  \sum_{F\in\mathbb F_G}
  \sum_{k=0}^p \gamma_k \left(\frac{h_F^{2k+1}}{(k!)^2} \llbracket \partial_n^k u_\ell \rrbracket, \llbracket \partial_n^k v_\ell \rrbracket\right)_F.
\end{equation*}
Here, \(\llbracket \cdot \rrbracket\) denotes the jump across face \(F\), \(p\) is the polynomial degree of the finite element space, \(h_F\) is the local face size, and \(\gamma_k\) are the ghost penalty parameters. Hence, given a right hand side function $f\in L^2(\Omega)$, our goal is to solve the weakly posed problem: find $u\in V_L$ such that
\begin{gather}
    \label{eq:weak-equation}
    A_L(u_L, v_L) = \left(f,v_L\right)_{\Omega},
    \qquad \forall v_L\in V_L.
\end{gather}
Below, we will also use the notation as a linear system of the form
\begin{gather}
    \label{eq:matrix-equation}
    \vA_L \vx_L = \vb_L.
\end{gather}
For $0\le \ell \le L$ we denote by $\vx_\ell$ the coefficient vector of the finite element function $u_\ell$ with respect to the finite element basis $\{\phi_{\ell,i}\}_{i=1,\dots,n_\ell}$, $\vb_\ell$ is the vector with entries $\left(f,\phi_{\ell,i}\right)$, and $\vA_\ell$ is the matrix with entries $A_\ell(\phi_{\ell,j},\phi_{\ell,i})$, where $i,j=1,\dots,n_\ell$.
We point out that only basis functions associated to the cells in $\Omega_\ell$ are part of this basis and no degrees of freedom are introduced on exterior cells.

\subsection{Multigrid}
We solve the linear system~\eqref{eq:weak-equation} either directly by a V-cycle iteration or using the V-cycle as a preconditioner for a GMRES method. The V-cycle is described in detail in~\cite{Bramble93,Hackbusch85}. It combines a so-called smoother on each level $\ell=1,\dots,L$ to be discussed in the following subsection with a recursion to the next coarser level $\ell-1$ and an exact solver on the coarsest level $\ell=0$.

In order to implement the recursion, we have to define a prolongation operator $\prol{\ell}$ transfering finite element functions from mesh $\mesh_{\ell-1}$ to mesh $\mesh_\ell$ and a restriction operator $\rest{\ell-1}$ transfering the residual from mesh $\mesh_{\ell}$ to mesh $\mesh_{\ell-1}$.

For fitted, nested meshes, it is customary to use the embedding operator from $V_{\ell-1}$ to $V_\ell$ as prolongation operator. In order to do so for unfitted meshes, we distinguish two cases, see also~\cite{LudescherGrossReusken20}. First the case of exact representation of the boundary: in this case, if a cell of $\mesh_\ell$ is cut by $\partial\Omega$, so will be coarse grid cell which contains it. Hence, it is easy to see that there holds $\Omega_{\ell} \subseteq \Omega_{\ell-1}$.
The second case arises if the boundary $\partial\Omega$ is approximated on every level, for instance by a finite element level set function. In that case, we follow~\cite{LudescherGrossReusken20} and introduce as an assumption what we had concluded in the first case, namely that $\Omega_{\ell} \subseteq \Omega_{\ell-1}$.
Hence, we can summarize that for any finite element function $u_{\ell-1}\in V_{\ell-1}$, the prolongation operator is well-defined by the simple condition
\begin{gather}
    \prol{\ell}u_{\ell-1}(x) = u_{\ell-1}(x) \quad\forall x\in \Omega_\ell.
\end{gather}
In our experiments, this assumption holds. For very general approximations, there might always be exceptions. At least, it would be an assumption of saturation type holding only for sufficiently fine meshes, which is not sufficient for multigrid methods.
If on the other hand the coarse grid level set is obtained by a suitable interpolation, it might be possible to prove this assumption. This is subject to further investigation.

The restriction operator is usually obtained as the adjoint of the prolongation operator with respect to a suitable inner product. The inner product we choose on $V_\ell$ is the Euclidean inner product of the coefficient vectors with respect to the nodal basis. Hence, in matrix representation, the restriction is the transpose of the prolongation.

\subsection{Vertex-patch smoother}

The vertex-patch smoother is a domain decomposition method employing an overlapping covering of $\Omega$, here actually $\Omega_\ell$, by very simple subdomains $\Omega_{\ell,j}$.

Let $\vertices_\ell$ be the set of vertices of $\mesh_\ell$ contained in $\Omega$. Let $\vertex_j\in\vertices_\ell$ be such a vertex. Then, the vertex patch $\mesh_{\ell,j}$ consists of all mesh cells in $\mesh_\ell$ having $\vertex_j$ as one of their vertices. Let $\Omega_{\ell,j}$ be the subdomain obtained as union of these cells. Note that by this very definition
\begin{gather*}
    \forall
    \vertex_j\in\vertices_\ell
    \quad\forall \cell \in \mesh_{\ell,j} 
    \quad\colon\quad \cell \subset \Omega_\ell.
\end{gather*}
Hence for all vertices $\vertex_j\in\vertices_\ell$ there holds $\Omega_{\ell,j} \subset\Omega_\ell$, and consequently, functions in the space $V_\ell$ are defined on the whole patch. This allows us to associate subspaces with each patch.

An extended version of the following discussion can be found in~\cite{CuiGrosseBleyKanschatStrzodka25,WichrowskiMunchKronbichlerKanschat25}.
For continuous finite element methods we define the subspaces $V_{\ell,j} \subset V_\ell$ as those functions with support on $\Omega_{\ell,j}$.
This implies in particular that the functions in $V_{\ell,j}$ vanish on the boundary of $\Omega_{\ell,j}$.
We refer to these as internal degrees of freedom of the patch.
The integrals in the bilinear form $a_\ell(\cdot,\cdot)$ in~\eqref{eq:cutfem_nitsche} couple these functions with themselves and those which are nonzero on the boundary of $\Omega_{\ell,j}$. We call those the boundary degrees of freedom. Internal and boundary degrees of freedom together form a larger subspace $\overline V_{\ell,j}$. Both subspaces are obtained by selecting a certain number of basis functions from the set $\{\phi_{\ell,i}\}$, hence the projection matrices $\vQ_{\ell,j}$ and $\overline \vQ_{\ell,j}$ from coefficient vectors in $V_\ell$ to those of $V_{\ell,j}$ and $\overline V_{\ell,j}$, respectively, simply select entries from a vector $\vx_\ell$. The projected matrices $\vA_{\ell,j}$ and $\overline \vA_{\ell,j}$ are then obtained by selecting the corresponding rows and columns from the matrix $\vA_\ell$.
We summarize this paragraph by the following representation of the so-called local solver
\begin{gather}
    \label{eq:local-solver}
    \vQ_{\ell,j}^T \vA_{\ell,j}^{-1}\vQ_{\ell,j} 
    \bigl(\vA_\ell\vx_\ell - \vb_\ell\bigr)
    = \vQ_{\ell,j}^T \vA_{\ell,j}^{-1}\vQ_{\ell,j} 
    \bigl(\overline \vA_{\ell_j} \overline\vQ_{\ell,j}\vx_\ell - \overline \vQ_{\ell,j}\vb_\ell\bigr).
\end{gather}
A single step of the multiplicative vertex-patch smoother (MVS) denoted by $S(\vx,\vb)$ and applied to a vector $\vx_\ell$ and right hand side $\vb_\ell$ consists of the steps
\begin{gather}
  \label{eq:smoother}
  \begin{split}
    \vy_\ell^{(0)} &= \vx_\ell,
    \\
    \vy_\ell^{(j)} &= \vQ_{\ell,j}^T \vA_{\ell,j}^{-1}\vQ_{\ell,j} 
    \Bigl(\vA_\ell\vy_\ell^{(j-1)} - \vb_\ell\Bigr),
    \qquad j=1,\dots,J,
    \\
    S(\vx,\vb) &= \vy_\ell^{(J)},   
  \end{split}
\end{gather}
where $J$ is the number of vertex patches.
Formula~\eqref{eq:local-solver} has two implications: first, the computation of the residual $\vA_\ell\vx_\ell - \vb_\ell$ and the solution of the local problem $\vA_{\ell,j}^{-1}$ can be obtained with a single sweep through the vector, reducing the memory footprint of the algorithm by a factor 2. Second, the results of two local solvers are independent of each other as soon as their patches are nonoverlapping. Hence, we can use a ``coloring'' of the patches into sets of nonoverlapping patches to execute local solvers in parallel. See~\cite{WitteArndtKanschat21,WichrowskiMunchKronbichlerKanschat25} for more details of this algorithm and on its cache coherence.

The ghost penalties in equation~\eqref{eq:cutfem-ghost} introduce a coupling between cells beyond degrees of freedom on the interface. Hence, the local solvers on two neighboring patches are not anymore independent. We have discussed such cases in~\cite{CuiKanschat24IP}, where we introduced an inconsistent version of the parallelized smoother and presented experiments, which showed that this version can be quite effective. Therefore, we decided to go the same way here and ignore coupling through ghost penalties between neighboring vertex patches.

\section{Implementation Details}

The finite element spaces are constructed using the so-called $hp$-framework of the deal.II finite element library~\cite{dealii} following step 85 in the tutorial, see also~\cite{StickoKreiss16}. Hence, degrees of freedom are only allocated on the cells in the level meshes $\mesh_{\ell,\Omega}$. The basis functions are those described after equation~\eqref{eq:spaces}.

Two main components of multigrid solvers are the operator action \(\vy \gets \vA\vx\) and the smoothing operation \(\vy \gets S(\vx,\vb)\), where the former represents the matrix-vector product in abstract form. In matrix-free algorithms, both the operator evaluation and the smoothing step are implemented via loops over cells or patches. In the following, we distinguish between interior patches, where all cells are subsets of $\Omega_\ell$, and cut patches, where at least one cell is cut by $\partial\Omega$. Interior patches have a Cartesian tensor product structure, cut patches do not.

Operator applications are implemented on the fly, using the sum factorization in~\cite{KronbichlerKormann12,KronbichlerKormann19} in the patch-wise pattern of~\cite{CuiKanschat24IP} on interior patches.
On cut patches, the tensor product structure is lost. Here, we have to use a substitute quadrature of sufficient accuracy, which at least on cut cells is not in tensor-product form. We employ the algorithm from~\cite{Saye2015} to generate these quadrature rules on the intersected cells. These computations are done upfront on the CPU.

Similarly, local solvers benefit from the Cartesian structure, since we can apply the fast diagonalization method like in~\cite{GottliebOrszag77,WitteArndtKanschat21,CuiGrosseBleyKanschatStrzodka25}, yielding a local solver with only $\mathcal O(N^{d+1})$ operations in $\mathbb R^d$, where $N$ is the number of degrees of freedom in direction of a one-dimensional tensor fiber, typically about twice the polynomial degree.
However, when a patch is cut, the inherent tensor structure is lost, rendering fast diagonalization infeasible. In such cases, a direct solver for the linear system of dimension $N^d$ in~\eqref{eq:local-solver} is employed. We compute the local matrices and their singular value decompositions upfront and load them to the GPU. 

Hence, the operator application and the smoother consist of ``structured'' and ``unstructured'' components. We provide details for the smoother~\eqref{eq:smoother}. The index set $1,\dots,J$ of the vertices is split into the interior and cut subsets $\mathbb J_i$ and $\mathbb J_c$, respectively, yielding the partitioned smoother given $\vx_\ell$:
 \begin{gather}
  \label{eq:smoother-split}
  \begin{split}
    \vx_\ell &\gets \vQ_{\ell,j}^T \vA_{\ell,j}^{-1}\vQ_{\ell,j} 
    \Bigl(\vA_\ell\vx_\ell - \vb_\ell\Bigr),
    \qquad j\in \mathbb J_i,
    \\
    \text{repeat }n_c\text{ times: }
    \vx_\ell &\gets \vQ_{\ell,j}^T \vA_{\ell,j}^{-1}\vQ_{\ell,j} 
    \Bigl(\vA_\ell\vx_\ell - \vb_\ell\Bigr),
    \qquad j\in \mathbb J_c,
    \\
    S(\vx_\ell,\vb_\ell) &\gets \vx_\ell.
  \end{split}
\end{gather}
In this algorithm, we already introduced the option of repeating the smoothing step on the cut cells, inspired by the domain decomposition method in~\cite{GrossReusken23} and the smoother in~\cite{LudescherGrossReusken20}.
Again, traversing the two subsets is parallelized by coloring.
We apply the highly optimized implementations from~\cite{CuiGrosseBleyKanschatStrzodka25,Cui24,CuiKanschat24IP} to the interior vertex patches in the subset $\mathbb J_i$, boosting performance of the overall algorithm. On the set $\mathbb J_c$ of cut patches, we employ a much less optimized algorithm. Under the assumption, that the size of the set $\mathbb J_i$ grows faster under mesh refinement than the size of the set $\mathbb J_c$ of patches at the boundary, we expect that the slowdown due to the less optimized code will eventually be mitigated.
\section{Computational Results}

In this section, we present numerical evidence demonstrating that the smoothers discussed above yield highly efficient multigrid methods, as reflected in low and stable iteration counts and high utilization of the hardware. All results are computed with the following setup: the coarsest mesh ($\ell=0$) consists of a 2-by-2 subdivision of the square $[-1.21,1.21]^2$ and all further meshes are obtained such that each cell of $\mesh_{\ell-1}$ is the union of four squares of $\mesh_\ell$. The domain $\Omega$ is the circle of radius $r=1$ around the origin. We use an analytic level set function with exact zero level on the perimeter of the circle. We run the standard V-cycle with a single pre- and post-smoothing step, respectively.

\begin{table}[tp]
\caption{Comparison of the vertex-patch smoother (MVS) on different geometries in 2D (left two blocks); comparison of MVS and a Chebyshev smoother (right two blocks). Multigrid preconditioner for GMRES solver, iteration counts to reduce residual by $10^{-9}$.}
\centering
\begin{tabular}{cccc|ccc|ccc}
\toprule
& \multicolumn{3}{c}{Square (MVS)} & \multicolumn{3}{c}{Circle (MVS)} & \multicolumn{3}{c}{Circle (Cheb.)} \\
  \cmidrule(lr){2-4} \cmidrule(lr){5-7} \cmidrule(lr){8-10}
 $L$ & $\mathbb{Q}_1$ & $\mathbb{Q}_2$ & $\mathbb{Q}_3$ & $\mathbb{Q}_1$ & $\mathbb{Q}_2$ & $\mathbb{Q}_3$ & $\mathbb{Q}_1$ & $\mathbb{Q}_2$ & $\mathbb{Q}_3$ \\
\midrule
6 & 5 & 4 & 3 & 6 & 9 & 17 & 4 & 11 & 65 \\
7 & 5 & 4 & 3 & 6 & 8 & 14 & 4 & 13 & 65 \\
8 & 5 & 4 & 3 & 6 & 7 & 13 & 4 & 12 & 59 \\
9 & 5 & 4 & 3 & 6 & 7 & 13 & 4 & 12 & 53 \\
10 & 5 & 4 & 3 & 5 & 7 & 14 & 4 & 12 & 53 \\
11 & 5 & 4 & 3 & 5 & 6 & 12 & 4 & 11 & 47 \\
12 & 5 & 4 & 3 & 5 & 5 & 10 & 4 & 12 & 46 \\
\bottomrule
\end{tabular}
 \label{tab:cut_conv}
\end{table}
First, we discuss the numerical efficiency, measured in iteration counts.
In Table~\ref{tab:cut_conv} on the left, we compare iteration counts between the multiplicative vertex-patch smoother (MVS) on fitted Cartesian meshes on the unit square as our baseline and the method described above on a circle. We see that we loose a few iterations by introducing cut cells, but the method converges reliably in very few steps.
On the domain with cut cells, we used the smoother with $n_c=2$ repetitions on the cut cells. 
On the right, we compare the MVS to a Chebyshev smoother  using Chebyshev polynomials of degree 5 for pre- and postsmoothing.
Here, the vertex-patch smoother can play out its advantage for higher order elements, as it is more robust.

Next, we discuss the influence of additional sweeps over the cut patches to the smoother.
This is motivated by the domain decomposition approach in~\cite{GrossReusken23} and also by the multigrid method in~\cite{LudescherGrossReusken20}. The results for the V-cycle preconditioner of the GMRES method are in Table~\ref{tab:cut_conv_cmp_gmres}.
\begin{table}[tp]
\caption{Vertex-patch smoother (MVS) with different number of smoothing steps $n_c$ for cut patches on circle domain in 2D. Multigrid preconditioner for GMRES solver, iteration counts to reduce residual by $10^{-9}$.}
\centering
\begin{tabular}{ccccccccccccc}
\toprule
 $n_c$ & 1 & 2 & 3 & 4 & 1 & 2 & 3 & 4 & 1 & 2 & 3 & 4\\
  \cmidrule(lr){2-5} \cmidrule(lr){6-9} \cmidrule(lr){10-13}
$L$ & \multicolumn{4}{c}{$\mathbb{Q}_1$} & \multicolumn{4}{c}{$\mathbb{Q}_2$} & \multicolumn{4}{c}{$\mathbb{Q}_3$} \\
\midrule
6 & 8 & 6 & 6 & 6   & 11 & 9 & 9 & 9    & 132 & 17 & 17 & 17 \\
7 & 8 & 6 & 6 & 6   & 9 & 8 & 8 & 8     & 193 & 14 & 14 & 14 \\
8 & 7 & 6 & 6 & 6   & 8 & 7 & 7 & 7     & 217 & 13 & 13 & 13 \\
9 & 7 & 6 & 5 & 5   & 8 & 7 & 7 & 7     & 402 & 13 & 13 & 13 \\
10 & 7 & 5 & 5 & 5  & 8 & 7 & 7 & 7     & 361 & 14 & 13 & 13 \\
11 & 7 & 5 & 5 & 5  & 8 & 6 & 6 & 6     & 307 & 12 & 12 & 13 \\
12 & 6 & 5 & 5 & 5  & 7 & 5 & 5 & 5     & 249 & 10 & 10 & 10 \\
\bottomrule
\end{tabular}
 \label{tab:cut_conv_cmp_gmres}
\end{table}
They indicate that a second iteration over all cut cells indeed reduces the number of iterations, albeit moderately for polynomial degrees one and two, such that in these cases the additional work does not pay off, unless on very fine meshes, where the work on the boundary vanishes compared to the bulk. The table also shows that a third iteration over the cut cells does not improve the convergence noticeably.
For $\mathbb Q_3$, the situation is different, and a second smoothing step on the cut cells is crucial for obtaining reasonable convergence rates. Unfortunately, our current implementation does not allow for higher polynomial degrees, as it relies on the computation of derivatives in deal.II~\cite{dealii}.
This discrepancy is even more pronounced when multigrid is applied as an iterative smoother itself and not as a preconditioner in GMRES.
Table~\ref{tab:cut_conv_cmp_vcycle} shows slightly higher iteration counts for $\mathbb Q_1$ and roughly doubled counts for $\mathbb Q_2$. For $\mathbb Q_3$ we observe divergence if $n_c=1$ and very high counts for more smoothing steps on the cut cells.
\begin{table}[tp]
\caption{Number of iterations for  vertex-patch smoother (MVS) with different number of smoothing steps $n_c$ for cut patches on the circle domain in 2D . Multigrid V-cycle solver with relative accuracy of $10^{-9}$.}
\centering
\begin{tabular}{ccccccccccccc}
\toprule
 $n_c$ & 1 & 2 & 3 & 4 & 1 & 2 & 3 & 4 & 1 & 2 & 3 & 4\\
  \cmidrule(lr){2-5} \cmidrule(lr){6-9} \cmidrule(lr){10-13}
$L$ & \multicolumn{4}{c}{$\mathbb{Q}_1$} & \multicolumn{4}{c}{$\mathbb{Q}_2$} & \multicolumn{4}{c}{$\mathbb{Q}_3$} \\
\midrule
6 & 14 & 10 & 8 & 8   & 33 & 24 & 23 & 23    & --- & 121 & 119 & 118 \\
7 & 14 & 9 & 8 & 8   & 22 & 18 & 17 & 17     & --- & 75 & 73 & 73 \\
8 & 12 & 8 & 8 & 8   & 18 & 13 & 13 & 12     & --- & 47 & 47 & 47 \\
9 & 12 & 8 & 8 & 7   & 18 & 13 & 13 & 13     & --- & 59 & 58 & 58 \\
10 & 12 & 8 & 8 & 7  & 22 & 14 & 13 & 13     & --- & 79 & 77 & 77 \\
11 & 12 & 8 & 7 & 7  & 22 & 14 & 13 & 13     & --- & 73 & 72 & 71 \\
12 & 12 & 8 & 7 & 7  & 22 & 13 & 13 & 13     & --- & 72 & 71 & 70 \\
\bottomrule
\end{tabular}
 \label{tab:cut_conv_cmp_vcycle}
\end{table}
 This indicates that while the spectral properties of the preconditioned operator are sufficient for GMRES, its largest or smallest eigenvalue prevents the linear method from converging fast.
We currently have no explanation for this behavior.

Finally, in the three tables~\ref{tab:cut_conv_ghost_parameter_q1} to~\ref{tab:cut_conv_ghost_parameter_q3}, we present a parameter study of the ghost penalty. 
\begin{table}[tp]
\caption{Number of iterations for  vertex-patch smoother (MVS) with various ghost penalty parameters $\gamma_1$ for $\mathbb{Q}_1$ element on the circle domain in 2D . Multigrid preconditioner for GMRES solver, iteration counts to reduce residual by $10^{-9}$.}
\centering
\begin{tabular}{c|cccccccccccc}
\toprule
 $L$ & 0.05 & 0.06 & 0.07 & 0.08 & 0.09 & 0.10 & 0.11 & 0.12 & 0.13 & 0.14 & 0.15 \\
\midrule
6 &  5.4 & 5.5 & 5.6 & 5.6 & 5.7 & 5.7 & 5.7 & 5.8 & 5.8 & 5.8 & 5.8 \\
7 &  5.0 & 5.2 & 5.2 & 5.3 & 5.3 & 5.3 & 5.4 & 5.4 & 5.4 & 5.5 & 5.5 \\
8 &  4.9 & 4.9 & 5.0 & 5.1 & 5.2 & 5.2 & 5.3 & 5.3 & 5.3 & 5.4 & 5.4 \\
9 &  4.9 & 4.8 & 4.8 & 4.9 & 4.9 & 4.9 & 5.0 & 5.2 & 5.2 & 5.3 & 5.3 \\
10 &  4.7 & 4.7 & 4.8 & 4.8 & 4.9 & 4.9 & 4.9 & 5.0 & 5.2 & 5.2 & 5.2 \\
11 &  7.3 & 6.5 & 5.8 & 4.9 & 4.8 & 4.9 & 4.9 & 5.0 & 5.1 & 5.2 & 5.2 \\
12 &  7.4 & 5.8 & 4.8 & 4.8 & 4.8 & 4.9 & 4.9 & 5.0 & 5.1 & 5.2 & 5.2 \\
\bottomrule
\end{tabular}
 \label{tab:cut_conv_ghost_parameter_q1}
\end{table}
\begin{table}[tp]
\caption{Number of iterations for  vertex-patch smoother (MVS) with various ghost penalty parameters $\gamma_2$ for $\mathbb{Q}_2$ element on the circle domain in 2D . Multigrid preconditioner for GMRES solver, iteration counts to reduce residual by $10^{-9}$.}
\centering
\begin{tabular}{c|cccccccccccc}
\toprule
 $L$ & 0.05 & 0.06 & 0.07 & 0.08 & 0.09 & 0.10 & 0.11 & 0.12 & 0.13 & 0.14 & 0.15 \\
\midrule
6 & 8.7 & 8.3 & 8.6 & 8.7 & 8.6 & 8.6 & 8.6 & 8.7 & 8.8 & 8.9 & 9.0 \\
7 & 7.3 & 7.3 & 7.4 & 7.4 & 7.5 & 7.6 & 7.7 & 7.8 & 7.8 & 7.8 & 7.8 \\ 
8 & 6.5 & 6.6 & 6.6 & 6.5 & 6.5 & 6.5 & 6.6 & 6.6 & 6.7 & 6.7 & 6.8 \\
9 & 5.9 & 6.0 & 6.2 & 6.2 & 6.3 & 6.3 & 6.4 & 6.4 & 6.5 & 6.5 & 6.5 \\
10 & 6.4 & 6.4 & 6.4 & 6.4 & 6.4 & 6.4 & 6.4 & 6.5 & 6.5 & 6.5 & 6.6 \\
11 & 5.6 & 5.6 & 5.6 & 5.6 & 5.7 & 5.7 & 5.7 & 5.7 & 5.7 & 5.7 & 5.7 \\
12 & 4.9 & 4.9 & 4.9 & 4.9 & 4.9 & 4.9 & 4.9 & 4.9 & 4.9 & 4.9 & 4.9 \\
\bottomrule
\end{tabular}
 \label{tab:cut_conv_ghost_parameter_q2}
\end{table}
\begin{table}[tp]
\caption{Number of iterations for  vertex-patch smoother (MVS) with various ghost penalty parameters $\gamma_3$ for $\mathbb{Q}_3$ element on the circle domain in 2D . Multigrid preconditioner for GMRES solver, iteration counts to reduce residual by $10^{-9}$.}
\centering
\begin{tabular}{c|cccccccccccc}
\toprule
 $L$ & 0.05 & 0.06 & 0.07 & 0.08 & 0.09 & 0.10 & 0.11 & 0.12 & 0.13 & 0.14 & 0.15 \\
\midrule
6 & 16.5 & 16.2 & 16.3 & 15.9 & 16.3 & 16.5 & 16.5 & 16.6 & 16.6 & 16.6 & 16.6 \\
7 & 13.9 & 13.8 & 13.6 & 12.8 & 13.6 & 13.8 & 13.8 & 13.8 & 13.8 & 13.8 & 13.8 \\
8 & 12.7 & 12.6 & 12.4 & 12.3 & 12.3 & 12.3 & 12.2 & 12.2 & 12.0 & 12.0 & 11.9 \\
9 & 13.7 & 13.5 & 13.0 & 12.9 & 12.8 & 12.7 & 12.6 & 12.6 & 12.5 & 12.5 & 12.5 \\
10 & 14.6 & 13.9 & 13.8 & 13.7 & 13.6 & 13.5 & 13.5 & 13.0 & 12.9 & 12.9 & 12.9 \\
11 & 12.6 & 12.0 & 11.9 & 11.9 & 11.8 & 11.8 & 11.7 & 11.7 & 11.7 & 11.6 & 11.6 \\
12 & 10.0 & 9.9 & 9.9 & 9.9 & 9.9 & 9.8 & 9.8 & 9.8 & 9.8 & 9.7 & 9.7 \\
\bottomrule
\end{tabular}
 \label{tab:cut_conv_ghost_parameter_q3}
\end{table}
For more fine-grained analysis, we use the fractional iteration counts $n_{\text{frac}}$ computed from the actual iterations $n_{\text{it}}$ by
\begin{gather*}
    n_{\text{frac}} = n_{\text{it}} \frac{-8} {\log_{10} \left(\|r_{n_{\text{it}}}\|/\|r_0\| \right)}.
\end{gather*}
For $\mathbb{Q}_1$, we search $\gamma_1$ over the interval $[0.05,0.15]$ and select the optimal value. For $\mathbb{Q}_2$, we fix $\gamma_1$ at the $\mathbb{Q}_1$ optimum and vary $\gamma_2$ to identify the best setting; $\mathbb{Q}_3$ follows the same procedure. We find that $\gamma_1$ exerts a stronger overall influence: if $\gamma_1$ is set suboptimally, no choice made in $\mathbb{Q}_2$ or $\mathbb{Q}_3$ is able to surpass the current best performance.

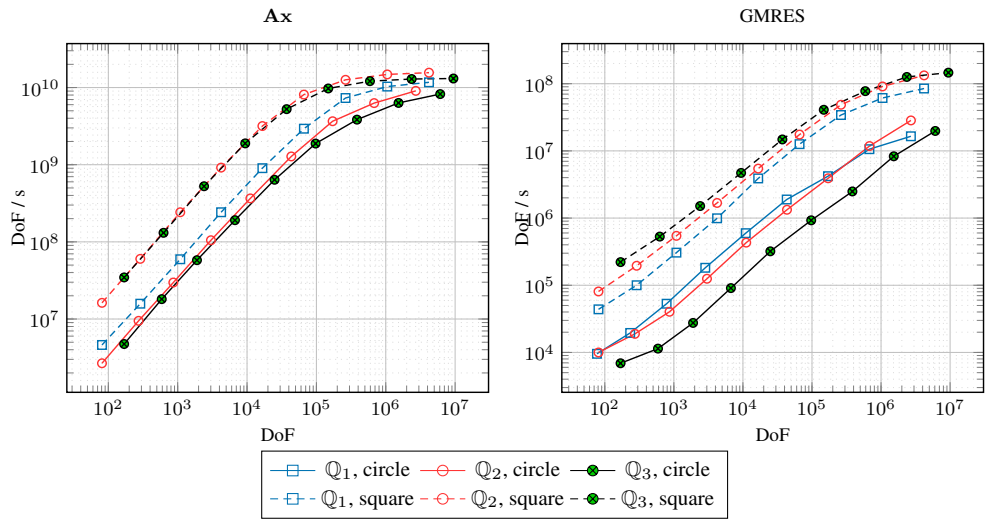
\begin{figure}[tp]
\centering
\footnotesize
\begin{tikzpicture}[scale=0.9]
\begin{axis}[
	xlabel={DoF},
	ylabel={DoF / s},
        ylabel shift=-2ex,
        ymode=log,
        log basis y={10},
        xmode=log,
        log basis x={10},
              grid=both,                        
      major grid style={gray!50},       
      minor grid style={gray!30,dotted},
      minor x tick num=8,               
      minor y tick num=8,
        width=0.6\textwidth,
        title={$\vA\vx$},
      legend columns=3,
      legend style={font=\footnotesize},
      cycle list name=colorGPL,
      mark size=1.8,
      semithick,
]
\addplot [draw=red!80!white,mark=o,every mark/.append style={solid}] table [x=v_q2,y=v_q2_p] {\cutfem}; 
\addplot [draw=black,mark=otimes*,every mark/.append style={solid,fill=green!80!black}] table [x=v_q3,y=v_q3_p] {\cutfem};

\addplot [draw=gnuplot@darkblue,densely dashed,mark=square,every mark/.append style={solid}] table [x=q1,y=q1_p] {\cutfem}; 
\addplot [draw=red!80!white,densely dashed,mark=o,every mark/.append style={solid}] table [x=q2,y=q2_p] {\cutfem};
\addplot [draw=black,mark=otimes*,densely dashed,every mark/.append style={solid,fill=green!80!black}] table [x=q3,y=q3_p] {\cutfem};
\end{axis}
\end{tikzpicture}
\hfill
\begin{tikzpicture}[scale=0.9]
\begin{axis}[
	xlabel={DoF},
	ylabel={DoF / s},
        ylabel shift=-2ex,
        ymode=log,
        log basis y={10},
        xmode=log,
        log basis x={10},
                      grid=both,                        
      major grid style={gray!50},       
      minor grid style={gray!30,dotted},
      minor x tick num=8,               
      minor y tick num=8,
        width=0.6\textwidth,
        title={GMRES},
        legend to name=legendCut,
      legend columns=3,
      legend style={font=\footnotesize},
      cycle list name=colorGPL,
      mark size=1.8,
      semithick,
]
\addplot [draw=gnuplot@darkblue,mark=square] table [x=dof_q1,y=gmres_q1] {\cutfemp}; 
\addlegendentry{$\mathbb{Q}_1$, circle}
\addplot [draw=red!80!white,mark=o,every mark/.append style={solid}] table [x=dof_q2,y=gmres_q2] {\cutfemp}; 
\addlegendentry{$\mathbb{Q}_2$, circle}
\addplot [draw=black,mark=otimes*,every mark/.append style={solid,fill=green!80!black}] table [x=dof_q3,y=gmres_q3] {\cutfemp};
\addlegendentry{$\mathbb{Q}_3$, circle}

\addplot [draw=gnuplot@darkblue,densely dashed,mark=square,every mark/.append style={solid}] table [x=dof_q1,y=gmres_q1] {\standardcutfemp}; \addlegendentry{$\mathbb{Q}_1$, square}
\addplot [draw=red!80!white,densely dashed,mark=o,every mark/.append style={solid}] table [x=dof_q2,y=gmres_q2] {\standardcutfemp};
\addlegendentry{$\mathbb{Q}_2$, square}
\addplot [draw=black,mark=otimes*,densely dashed,every mark/.append style={solid,fill=green!80!black}] table [x=dof_q3,y=gmres_q3] {\standardcutfemp};
\addlegendentry{$\mathbb{Q}_3$, square}
\end{axis}
\end{tikzpicture}
\ref{legendCut}
\caption{Throughput of matrix-free operator evaluation and solution on circle domain compared to the throughput on square domain.}
\label{fig:cutfem_perf}
\end{figure}

Finally, we assess the computational performance of the CutFEM method as the elapsed computation time on the Nvidia A100. Figure~\ref{fig:cutfem_perf} compares the throughput, measured in degrees of freedom per second (DoF/s), for matrix-free operator evaluation and for the whole GMRES solver on both square and circle domains. Note that on the same mesh, the number of degrees of freedom differs between the full square domain and the circle domain on the Cartesian background mesh, as the exterior cells do not carry degrees of freedom.
The left panel of the figure reveals a performance gap of about a factor seven on coarse meshes. This is due to the fact that the very fast quadrature by sum factorization is replaced by a more complicated procedure on cut cells.
However, as the mesh is refined beyond a million degrees of freedom and the proportion of cut cells decreases, the performance on the two domains becomes comparable. Furthermore, owing to our matrix-free implementation, the performance remains stable across different polynomial orders.

For the whole solver, however, a decrease in overall performance is observed on the right of Figure~\ref{fig:cutfem_perf} due to two primary factors. First, an increased iteration count directly lowers throughput. Hence, even for fine meshes, the gap between the fitted and unfitted meshes does not close. According to Table~\ref{tab:cut_conv}, this influence is small for the lowest order element, but already a factor of two for $\mathbb Q_2$.
Second, because the smoother employs a direct solver, which has a higher computational complexity than fast diagonalization and needs more data, the local matrices may be retrieved from L2 cache or even from global memory, where less efficient memory accesses can substantially impact performance.
Nevertheless, the performance gap is considerably less than a factor 10 for fine meshes, but it is about an order of magnitude on the coarsest meshes.
\pgfplotstableread{
l	v_q1	v_q1_p	v_q2	v_q2_p	v_q3	v_q3_p	q1	q1_p	q2	q2_p	q3	q3_p
4	77	9.13E+05	81	2.67E+06	169	4.74E+06	81	4.61E+06	81	1.62E+07	169	3.46E+07
5	233	2.63E+06	273	9.47E+06	589	1.81E+07	289	1.58E+07	289	6.02E+07	625	1.31E+08
6	785	8.10E+06	865	2.97E+07	1897	5.79E+07	1089	5.95E+07	1089	2.42E+08	2401	5.26E+08
7	2849	2.71E+07	3017	1.05E+08	6697	1.91E+08	4225	2.42E+08	4225	9.17E+08	9409	1.89E+09
8	10981	9.88E+07	11161	3.66E+08	24937	6.38E+08	16641	9.02E+08	16641	3.17E+09	37249	5.25E+09
9	43065	3.72E+08	43457	1.28E+09	97429	1.88E+09	66049	2.93E+09	66049	8.11E+09	148225	9.73E+09
10	170481	1.24E+09	171329	3.66E+09	384793	3.84E+09	263169	7.32E+09	263169	1.26E+10	591361	1.21E+10
11	678225	3.68E+09	680065	6.29E+09	1528753	6.32E+09	1050625	1.03E+10	1050625	1.48E+10	2362369	1.29E+10
12	2705341	6.08E+09	2709185	9.06E+09	6092881	8.23E+09	4198401	1.17E+10	4198401	1.56E+10	9443329	1.31E+10
}\cutfem
\pgfplotstableread{
dof_q1	it_q1	gmres_q1	dof_q2	it_q2	gmres_q2	dof_q3	it_q3	gmres_q3 h_1
77	7	9.54E+03	81	6	9.99E+03	169	9	6.90E+03 1.81
233	8	1.95E+04	273	8	1.89E+04	589	13	1.14E+04 3.62
785	8	5.31E+04	865	9	4.03E+04	1897	13	2.75E+04 7.24
2849	7	1.82E+05	3017	8	1.25E+05	6697	11	9.07E+04 14.47
10981	7	5.96E+05	11161	7	4.31E+05	24937	9	3.20E+05 28.95
43065	7	1.89E+06	43457	7	1.33E+06	97429	10	9.26E+05 57.8
170481	7	4.20E+06	171329	7	3.91E+06	384793	11	2.49E+06 116.28
678225	7	1.06E+07	680065	6	1.18E+07	1528753	9	8.31E+06 232.55
2705341	7	1.65E+07	2709185	5	2.84E+07	6092881	8	1.98E+07 454.54
}\cutfemp
\pgfplotstableread{
dof_q1	it_q1	gmres_q1	dof_q2	it_q2	gmres_q2	dof_q3	it_q3	gmres_q3 h_1
81	5	4.38E+04	81	4	8.08E+04	169	3	2.21E+05 1.81
289	6	1.00E+05	289	4	1.95E+05	625	3	5.30E+05 3.62
1089	6	3.06E+05	1089	4	5.44E+05	2401	3	1.51E+06 7.24
4225	6	9.92E+05	4225	4	1.68E+06	9409	3	4.70E+06 14.47
16641	5	3.89E+06	16641	4	5.45E+06	37249	3	1.48E+07 28.95
66049	5	1.26E+07	66049	4	1.75E+07	148225	3	4.09E+07 57.8
263169	5	3.42E+07	263169	4	4.86E+07	591361	3	7.75E+07 116.28
1050625	5	6.13E+07	1050625	4	9.10E+07	2362369	3	1.26E+08 232.55
4198401	5	8.48E+07	4198401	4	1.33E+08	9443329	3	1.46E+08 454.54
}\standardcutfemp

\section{Conclusions}

We have presented a multigrid method with vertex-patch smoother for unfitted finite element methods, in particular cutFEM with face-based ghost penalties. We discussed in particular the implementation on GPU, where we use code optimized for tensor product elements for the uncut vertex patches in the interior of the domain. As we loose tensor product structure on the cut patches, we employ a direct solver for the local problems there. Our results show, that the performance penalty is moderate in simple cases. We observe robustness under mesh refinement and the preconditioner also works fairly well for higher order.
As the next challenge we will have to address the inconsistent treatment of the ghost penalty in order to reduce iteration counts. In~\cite{CuiKanschat24IP}, we used Hermitian basis functions for this, which will not work here, since the ghost penalties involve high order derivatives. We will also investigate the suitability of different stabilization techniques in the context of vertex-patch smoothers.

The second challenge to reduce the performance gap will be faster solvers for the local problems on cut patches. We have good experience with iterative local solvers~\cite{CuiKanschat25Stokes}, see also~\cite{BastianMuellerMuethingPiatkowski19}, which is a direction to explore here as well.

\bibliographystyle{siam}
\bibliography{references}

\end{document}